\documentclass[12pt]{article}
\usepackage{latexsym,amsfonts,amssymb}
\usepackage[dvips]{color}
\usepackage{colordvi,multicol}

\newtheorem{thm}{Theorem}[section]

\newtheorem{lem}[thm]{Lemma}

\newtheorem{defi}[thm]{Definition}

\begin{document}

\title{\bf Chebyshev's inequality for
Banach-space-valued random elements}
\author{Ling Zhou and Ze-Chun Hu\thanks{Corresponding
author: Department of Mathematics, Nanjing University, Nanjing
210093, PR China\vskip 0cm E-mail address: huzc@nju.edu.cn}\\
 {\small Nanjing University}}
 \date{}
 \maketitle

%


\noindent{\bf Abstract}\quad In this paper, we obtain a new
generalization of Chebyshev's inequality for  random elements taking
values in a separate Banach space.

\smallskip


\smallskip



\section{Introduction}
Chebyshev's inequality states that for a random variable $X$ with
mean $E(X)$ and variance Var$(X)$, and any $\varepsilon>0$,
\begin{eqnarray*}
P\{|X-E(X)|\geq \varepsilon\}\leq \frac{{\rm
Var}(X)}{\varepsilon^2}.
\end{eqnarray*}
This inequality plays an important role in probability theory and
statistics. Several generalization for random vectors have been
made. A natural one is as follows.

Suppose that $X$ is an $n$-dimensional random vector, then for
$\varepsilon>0$,
\begin{eqnarray*}
P\{\|X-E(X)\|\geq \varepsilon\}\leq \frac{{\rm
Var}(X)}{\varepsilon^2},
\end{eqnarray*}
where $\|\cdot\|$ denotes the Euclidean norm in $\mathbf{R}^n$ and
Var$(X)=E[\|X-E(X)\|^2]$. This result can be seen in Laha and
Rohatgi (1979), P. 446-451. Other generalizations are given in
Marshall and Olkin (1960), Godwin (1955) and Mallows (1956).

Grenander (1963) proved a Chebyshev's inequality for
Hilbert-space-valued random elements as follows: if  $X$ is a random
element taking values in a Hilbert space $H$ with
$E(\|X\|^2)<\infty$, then for $\varepsilon>0$,
\begin{eqnarray*}
P\{\|X\|\geq \varepsilon\}\leq \frac{E(\|X\|^2)}{\varepsilon^2}.
\end{eqnarray*}

Chen (2007) proved the following new generalization of Chebyshev's
inequality for random vectors.

\begin{thm}\label{thm1.1}
Suppose that $X$ is an $n$-dimensional random vector with positive
definite covariance matrix $\Sigma$. Then, for any $\varepsilon>0$,
\begin{eqnarray*}
P\{(X-E(X))^{T}\Sigma^{-1}(X-E(X))\geq \varepsilon\}\leq
\frac{n}{\varepsilon},
\end{eqnarray*}
where the superscript ``\,$T$" denotes the transpose of a matrix.
\end{thm}

Rao (2010) extended Theorem \ref{thm1.1} to   random elements taking
values in a separable Hilbert space as follows.
\begin{thm}\label{thm1.2}
Suppose that $X$ is a random element taking values in a separable
Hilbert space H with expectation zero, positive definite covariance
operator $S$, and probability distribution $\mu$ such that $\int_H
{||x|{|^2}\mu(dx) }<\infty.$ Then, for every $\varepsilon >0$,
$$P\{(SX,X) > \varepsilon \} \le \frac{{{{\left[ \int_H {||x|{|^2}\mu(dx) } \right]}^2}}}{\varepsilon }$$
and
$$P\{(S^{-1}X,X) > \varepsilon \}
 \le \frac{{{{\left[ ||S^{-1}||\int_H {||x|{|^2}\mu(dx)} \right]}^2}}}{\varepsilon },$$
where covariance operator $S$ is the  Hermitian operator determined
uniquely by the quadratic form $(Sy,y) = \int_H {{{(x,y)}^2}}
\mu(dx).$
\end{thm}

In this paper, we will extend Theorem \ref{thm1.2} to random
elements taking values in a separable Banach space.

\section{Main result}
Suppose that $B$ is a real separable Banach space, $\|\cdot\|$ is
the norm, and $\mathcal {B}(B)$ is the Borel $\sigma$-algebra. $B^*$
is the dual space of $B$, i.e. $B^*$ is the family of all bounded
linear functionals on $B$. It's well known that $B^*$ is also a
Banach space with the operator norm $\|\cdot\|^*$ defined by
$$
\|f\|^*=\sup_{\|x\|\leq 1}|f(x)|,\ \forall f\in B^*.
$$
Let $(\Omega,\mathcal{F},P)$ be a probability space. A mapping $X:
(\Omega,\mathcal{F})\to (B,\mathcal{B}(B))$ is called measurable if
for any $A\in \mathcal{B}(B)$, we have $X^{-1}(A)\in\mathcal{F}$,
where $X^{-1}(A)=\{\omega\in\Omega: X(\omega)\in A\}$. We call such
$X$ random element taking values in $B$. Let $\mu$ be the
probability distribution of $X$, i.e. $\mu$ is the probability
measure on $(B,\mathcal{B}(B))$ defined by
$$
\mu(A)=P(X^{-1}(A)),\ \forall A\in \mathcal{B}(B).
$$

\begin{defi}\label{defi2.1}
Suppose that a mapping $X: (\Omega,\mathcal{F})\to
(B,\mathcal{B}(B))$ is measurable. $X$ is called a step function if
it can be expressed by
$$
X=\sum_{i=1}^{\infty}x_iI_{A_i},
$$
where  $\forall i\in \mathbf{N},x_i\in B,A_i\in\mathcal{F}$ and
$A_i\cap A_j=\emptyset,\forall i\neq j$.
\end{defi}

By the separability of $B$, we can easily get

\begin{lem}\label{lem2.2}
Suppose that $X: (\Omega,\mathcal{F})\to (B,\mathcal{B}(B))$ is
measurable. Then there exists a sequence $\{X_n,n\geq 1\}$ of step
functions such that
$$
\lim_{n\to\infty}\sup_{\omega\in\Omega}\|X_n(\omega)-X(\omega)\|=0
$$
and $\forall \omega\in\Omega,\forall n\geq 1,\ \|X_n(\omega)\|\leq
2\|X(\omega)\|$.
\end{lem}

To state our main result, we should define the covariance operator
of random element $X$ taking values in the Banach space $B$. For any
$x\in B,f\in B^*$,  define $(x,f)=(f,x)=f(x)$.

\begin{thm}\label{thm2.3}
Suppose that $X$ is a random element  taking values in $B$ with
probability distribution  $\mu$.  If $\int_B\|x\|^2\mu(dx)<\infty$,
then the quadratic form
\begin{eqnarray}\label{thm2.3-a}
(Sf,g)=\int_B(x,f)(x,g)\mu(dx),\ \forall f,g\in B^*,
\end{eqnarray}
uniquely determines one bounded linear operator $S: B^*\to B$.
\end{thm}
{\bf Proof.} {\bf Existence:} Firstly, we prove that  $\forall f\in
B^*, ~\exists ~Sf\in B$ such that
\begin{eqnarray}\label{thm2.3-b}
(Sf,g)=\int_B f(x)g(x)\mu(dx),\ \forall g\in B^*.
\end{eqnarray}
By Lemma \ref{lem2.2}, there exists a sequence $\{X_n,n\geq 1\}$ of
step functions such that $X_n$ converges to $X$ uniformly and
$\|X_n(\omega)\|\leq 2|X(\omega)|,\forall \omega\in\Omega,\forall
n\geq 1$. Let $\mu_n$ be the probability distribution of $X_n$. Then
for any $n\geq 1$, we have
$$
\int_B\|x\|^2\mu_n(dx)=\int_{\Omega}\|X_n\\|^2dP\leq
4\int_{\Omega}\|X\|^2dP=4\int_B||x||^2\mu(dx)<\infty.
$$
Let $X_1$ have the expression that $
X_1=\sum_{i=1}^{\infty}x_iI_{A_i}$ as in Definition \ref{defi2.1}.
For $f\in B^*$, define
\begin{eqnarray}\label{thm2.3-c}
S_1f=\sum_{i=1}^{\infty}f(x_i)P(A_i)x_i.
\end{eqnarray}
Since
\begin{eqnarray}\label{thm2.3-c-1}
\sum_{i=1}^{\infty}\|f(x_i)P(A_i)x_i\|&=&
        \int_\Omega\|f(X_1)X_1\|dP\nonumber\\
        &\leq&\|f\|\int_\Omega\|X_1\|^2dP\nonumber\\
        &=&\|f\|\int_B\|x\|^2\mu_1(dx)<\infty,
\end{eqnarray}
 $S_1f$ is well defined and we have
\begin{eqnarray}\label{thm2.3-d}
S_1f= \int_\Omega f(X_1)X_1dP.
\end{eqnarray}
 By (\ref{thm2.3-c}) and (\ref{thm2.3-c-1}), we get that for any
$g\in B^*$,
\begin{eqnarray*}
(S_1f,g)&=&\sum_{i=1}^{\infty}f(x_i)P(A_i)g(x_i)\\
&=&\int_\Omega f(X_1)g(X_1)dP\\
&=& \int_Bf(x)g(x)\mu_1(dx).
\end{eqnarray*}
For $n=2,3,\ldots$, define $S_nf$ from $X_n$ similar to $S_1f$. In
particular, we have
\begin{eqnarray}\label{thm2.3-e}
S_nf= \int_\Omega f(X_n)X_ndP,
\end{eqnarray}
and for any $g\in B^*$,
\begin{eqnarray}\label{thm2.3-e-1}
(S_nf,g)&=&\int_\Omega f(X_n)g(X_n)dP= \int_Bf(x)g(x)\mu_n(dx).
\end{eqnarray}
For any $n,m=1,2,\ldots$, by (\ref{thm2.3-d}) and (\ref{thm2.3-e}), we have
\begin{eqnarray}\label{thm2.3-f}
\|S_nf-S_mf\|&=&\left\|\int_{\Omega}\left(f(X_n)X_n-
f(X_m)X_m\right)dP\right\|\nonumber\\
&\leq&\int_{\Omega}\left\|f(X_n)X_n-f(X_m)X_m\right\|dP\nonumber\\
&=&\int_\Omega\|f(X_n)(X_n-X_m)+X_m(f(X_n)-f(X_m))\|dP\nonumber\\
&\le&\|f\|^*\int_\Omega(\|X_n\|+\|X_m\|)(\|X_n-X_m\|)dP.
\end{eqnarray}
By the fact that  $\|X_n(\omega)\|\leq 2\|X(\omega)\|,\forall n\geq
1,\forall \omega\in\Omega$, we have
$$
\int_{\Omega}(\|X_n\|+\|X_m\|)(\|X_n-X_m\|)dP\leq 16\int_{\Omega}\|X\|^2dP=16\int_B\|x\|^2\mu(dx)<\infty.
$$
Then it follows from (\ref{thm2.3-f}) and the dominated convergence theorem that \begin{eqnarray*}
\|S_nf-S_mf\|\to
0\ \ \mbox{as}\ n,m\to\infty,
\end{eqnarray*}
which implies that  $\{S_nf,n\geq 1\}$ is a Cauchy sequence in $B$. Thus
 there exists a unique element denoted by $Sf$ such that $S_nf$ converges to $Sf$ in $B$.
Furthermore, by (\ref{thm2.3-e-1}) and integral transformation, we have
\begin{eqnarray*}
(Sf,g)=\int_\Omega
f(X)g(X)dP=\int_B
f(x)g(x)\mu(dx),\ \forall g\in B^*.
\end{eqnarray*}

Secondly,  we prove that $S$ is {\bf linear}. Suppose that
$f_1,f_2\in B^*$ and $a,b\in \mathbf{R}$. Then for any $g\in B^*$,
\begin{eqnarray*}
(S(af_1+bf_2),g)&=&\int_B(x,af_1+bf_2)(x,g)\mu(dx)\\
                &=&a\int_B(x,f_1)(x,g)\mu(dx)+b\int_B(x,f_2)(x,g)\mu(dx)\\
                &=&(aSf_1,g)+(bSf_2,g)\\
                &=&(aSf_1+bSf_2,g),
\end{eqnarray*}
which implies that $S(af_1+bf_2)=aSf_1+bSf_2$.

Thirdly,  we prove that $S$ is {\bf bounded}. By (\ref{thm2.3-e}),
for any $n\geq 2$, we have
\begin{eqnarray}\label{thm2.3-g}
\|S_n f\|\le \|f\|^*\int_\Omega \|X_n\|^2dP.
\end{eqnarray}
Notice that $S_nf$ converges to $Sf$ and $X_n$ converges uniformly to $X$. Then
 letting $n\to \infty$ in (\ref{thm2.3-g}), we get
\begin{eqnarray*}
\|S f\|\le \|f\|^*\int_\Omega \|X\|^2dP=\|f\|^*\int_B\|x\|^2\mu(dx),
\end{eqnarray*}
which implies that $S$ is a bounded operator from $B^*$ to $B$.

{\bf Uniqueness:} Suppose that  $S':B^*\to B$ is another bounded
linear operator satisfying that
$$
(S'f,g)=\int_B(x,f)(x,g)\mu(dx),\forall
f,g \in B^*.
$$
Then for any $f\in B^*$, we have
$$
(S'f-Sf,g)=0,\ \forall g\in B^*,
$$
which implies that $S'f-Sf=0$. Thus $S'=S$.  \hfill\fbox

\bigskip

Suppose that $P^*$ is a probability measure on
$(B^*,\mathcal{B}(B^*))$. Since $S$ is a bounded linear operator
from $B^*$ to $B$, we can check that $f\mapsto (Sf,f)$ is a
nonnegative continuous functional on $B^*$, and thus it is
measurable with respect to $\mathcal{B}(B^*)$. For any
$\varepsilon>0,$ define
$$D_\varepsilon=\{f\in B^*:(Sf,f)\ge\varepsilon\}.$$
Then
\begin{eqnarray*}
P^*( D_\varepsilon)&\le& \frac{1}{\varepsilon }\int_{D_\varepsilon}(Sf,f)P^*(df)\\
&\le&\frac{1}{\varepsilon }\int_{B^*}(Sf,f)P^*(df)\\
&=&\frac{1}{\varepsilon }\int_{B^*}\left(\int_B f^2(x)\mu(dx)\right)P^*(df)\\
&\le&\frac{1}{\varepsilon }\int_{B^*}\left(\int_{B}(\|f\|^*)^2\|x\|^2\mu(dx)\right)P^*(df)\\
&=&\frac{1}{\varepsilon }\int_{B^*}(\|f\|^*)^2P^*(df)\left(\int_B{\|x\|^2} \mu(dx)\right).\\
\end{eqnarray*}

We have known that $S$ is nonnegative definite, i.e. for any $f\in
B^*$, $(Sf,f)\geq 0$. Furthermore, if $S$ is positive definite in
the sense that $(Sf,f)=0$ implies that $f=0$,  then $S$ is
invertible. For any $y\in B$, we have
$$(SS^{-1}y,S^{-1}y)=\int_B(x,S^{-1}y)^2\mu(dx)\geq 0,$$
i.e.
$$
(y,S^{-1}y)=\int_B(x,S^{-1}y)^2\mu(dx)\geq 0.$$ Define
$$
D_\varepsilon^{'}=\{y\in
B:(S^{-1}y,y)\ge\varepsilon\}.
$$
Then
\begin{eqnarray*}
 P\{X\in D_\varepsilon'\} &=&\int_{ D_\varepsilon^{'}} \mu(dy)  \\
  &\le& \frac{1}{\varepsilon }\int_{D_\varepsilon ^{'}} {(y,{S^{-1}}y)\mu(dy) }  \\
  &\le& \frac{1}{\varepsilon }\int_B (y,S^{ - 1}y)\mu (dy)  \\
  &=& \frac{1}{\varepsilon }\int_B \left(\int_B (x,S^{ - 1}y)^2\mu (dx)\right)\mu (dy)  \\
  &\le& \frac{1}{\varepsilon}\int_B \left(\int_B  \|x\|^2(\|S^{-1}y\|^*)^2\mu (dx) \right)\mu (dy)  \\
  &\le& \frac{1}{\varepsilon }\int_B \left( \int_B  \|x\|^2\|S^{-1}\|^2\mu (dx) \right)\|y\|^2\mu (dy)  \\
  &=& \frac{1}{\varepsilon }\|S^{-1}\|^2\left[ \int_B \|x\|^2\mu (dx)
  \right]^2,
\end{eqnarray*}
where $\|S^{-1}\|$ is the operator norm of $S^{-1}:B\to B^*$. Hence
we have the following result.

\begin{thm}\label{thm2.4}
Suppose that  $X$ is a random element taking values in $B$ with probability distribution
$\mu$ satisfying that $\int_B\|x\|^2\mu(dx)<\infty$, $P^*$ is
a  probability measure on $(B^*,\mathcal{B}(B^*))$,
 $S:{B}^*\to {B}$ is the bounded linear operator defined by  Theorem
\rm{\ref{thm2.3}}. Then for any $\varepsilon>0$, we have
$$P^*\{f\in B^*:(Sf,f)\ge\varepsilon\}\le\frac{1}{\varepsilon }
\int_{B^*}(\|f\|^*)^2P^*(df)\left(\int_B{\|x\|^2} \mu(dx)\right)$$
and
$$P\{(S^{-1}X,X)\ge\varepsilon\}\le\frac{1}{\varepsilon }\|S^{ - 1}\|^2\left[\int_B \|x|^2\mu (dx)\right]^2.$$
\end{thm}

\section{Remarks}\setcounter{equation}{0}
In this section, we show that Theorem \ref{thm2.4} extends Theorem
\ref{thm1.2}. Let $H$ be a separable Hilbert space with inner
product $(\cdot,\cdot)_H$ and  $H^*$ be the dual space of $H$ with
norm $\|\cdot\|^*$. Let $X$ be a random element taking values in $H$
with probability distribution $\mu$ satisfying
$\int_H\|x\|^2\mu(dx)<\infty$.  Let $S_H$ be the covariance operator
 defined in Theorem \ref{thm1.2}, and $S$ be the bounded linear
operator from $H^*$ to $H$ defined by Theorem \ref{thm2.3}.

By Riesz representation theorem, there exists an isometry $T: H\to H^*$. In fact, for any $x\in H$,
$$
Tx(y)=(y,x)_H,\ \forall y\in H,
$$
and $\|Tx\|^*=\|x\|$. As in Section 2, for any $x\in H,f\in H^*$, we
define $(x,f)=(f,x)=f(x)$. Then we have $S T=S_H$. In fact, for any
$y\in H$, we have $STy\in H$ and
\begin{eqnarray*}
(STy,y)_H&=&Ty(STy)=(S Ty,Ty)\\
&=&\int_H((Ty(x))^2\mu(dx)\\
&=&\int_H(x,y)^2_H\mu(dx)=(S_Hy,y)_H.
\end{eqnarray*}

Define one probability measure $P^*$ on $(H^*,\mathcal{B}(H^*))$ by
$$
P^*(A)=\mu(T^{-1}(A)),\ \forall A\in\mathcal{B}(H^*).
$$
Then
\begin{eqnarray}\label{3.1}
P^*(D_{\varepsilon})&=&\mu(T^{-1}(D_{\varepsilon}))=P(X^{-1}(T^{-1}(D_{\varepsilon})))\nonumber\\
&=&P\{\omega\in\Omega:T
X(\omega)\in
D_\varepsilon\}\nonumber\\
&=&P\{\omega\in\Omega:(STX(\omega),TX(\omega))\geq \varepsilon\}.
\end{eqnarray}
By the fact that $\|Tx\|^*=\|x\|,\forall x\in H$, and integral
transformation, we have
\begin{eqnarray}\label{3.2}
\int_H\|x\|^2\mu(dx)=\int_H(\|Tx\|^*)^2\mu(dx)=\int_{H^*}(\|f\|^*)^2P^*(df).
\end{eqnarray}
By Theorem \ref{thm2.4}, (\ref{3.1}) and (\ref{3.2}), we obtain
\begin{eqnarray}\label{3.3}
P\{\omega\in\Omega:(STX(\omega),TX(\omega))\geq \varepsilon\} &\leq&
\frac{1}{\varepsilon }
\int_{H^*}(\|f\|^*)^2P^*(df)\left(\int_H{\|x\|^2}
\mu(dx)\right)\nonumber\\
&=&\frac{1}{\varepsilon } \left(\int_H{\|x\|^2}\mu(dx)\right)^2.
\end{eqnarray}
By $ST=S_H$ and the definition of $T$, we have
\begin{eqnarray}\label{3.4}
(STX(\omega),TX(\omega))=(S_HX,X)_H.
\end{eqnarray}
It follows from (\ref{3.3}) and (\ref{3.4}) that
\begin{eqnarray}\label{3.5}
P\{(S_HX,X)_H\geq \varepsilon\}\leq  \frac{1}{\varepsilon }
\left(\int_H{\|x\|^2}\mu(dx)\right)^2.
\end{eqnarray}

On the other hand, if $S_H$ is positive definite, then $S$ is positive definite by
$ST=S_H$ and the isometry of $T$. Note that $S^{-1}X\in H^*$ and
$$
(S^{-1}X,X)=S^{-1}X(X)=TT^{-1}S^{-1}X(X)=(X,T^{-1}S^{-1}X)_H=(S_H^{-1}X,X)_H.
$$
Then we  have
\begin{eqnarray}\label{3.6}
P\{X\in D_{\varepsilon}'\}&=&P\{(S^{-1}X,X)\geq
\varepsilon\}=P\{(S_H^{-1}X,X)_H\geq \varepsilon\}.
\end{eqnarray}
By $S_H^{-1}=T^{-1}S^{-1}$ and the isometry of $T$, we have
$\|S^{-1}\|=\|S_H^{-1}\|$. By Theorem \ref{thm2.4} and (\ref{3.6}),
we obtain that
\begin{eqnarray}\label{3.7}
P\{(S_H^{-1}X,X)\geq \varepsilon\}\leq  \frac{1}{\varepsilon }
\|S_H^{-1}\|^2\left[\int_H \|x|^2\mu (dx)\right]^2.
\end{eqnarray}
Inequalities (\ref{3.5}) and (\ref{3.7}) are just those two ones in
Theorem \ref{thm1.2}.

\bigskip

{ \noindent {\bf\large Acknowledgments} \vskip 0.1cm  \noindent
Research supported by  NNSFC (Grant No. 10801072) and the
Fundamental Research Funds for the Central Universities.}




\begin{thebibliography}{90}\addcontentsline{toc}{section}{$\mathbf{References}$}
{
\bibitem{Chen07} Chen, X., 2007. A new generalization of Chebyshev
inequality for random vectors. arXiv: 0707.0805v1[math.ST] 5 Jul
2007.

\bibitem{Go55}Godwin, H.J., 1955. On generalizations of Chebyshev
inequality. Journal of the American Statistical Association 50,
923-945.

\bibitem{Gr63}Grenander, U., 1963. Probabilities on Algebraic
Structures. Wiley, New York, Almqvist and Wiksell, Stockholm.

\bibitem{LR79}Laha, R.G., Rohatgi, V.K., 1979. Probability Therory.
Wiley, New York.

\bibitem{Ma56}Mallows, C.L., 1956. Generalizations of Chebyshev
inequalities. Journal of Royal Statistical Society, Series B 18,
139-171.

\bibitem{MO60}Marshall, A., Olkin, I., 1960. Multivariate
Chebyshev inequalities. Annals of Mathematical Statistics 31,
1001-1014.


\bibitem{Rao10} Rao, B.L.S.P., 2010. Chebyshev's inequality
for Hilbert-space-valued random elements. Statistics and Probability
Letters 80, 1039-1042.
    }
\end{thebibliography}
\end{document}